\documentclass[12 pt]{article}
\usepackage{latexsym, amssymb, amsmath, amsfonts, amsthm,graphicx}
\theoremstyle{plain}
		\newtheorem{theorem}{Theorem}[section]
		\newtheorem{lemma}{Lemma}[section]

\theoremstyle{definition}
		\newtheorem{definition}{Definition}[section]

\theoremstyle{remark}

\def\beq{ \begin{equation} }
\def\eeq{ \end{equation} }
\def\mn{\medskip\noindent}

\def\ep{\epsilon}

\def\square{\vcenter{\vbox{\hrule height .4pt
  \hbox{\vrule width .4pt height 5pt \kern 5pt
        \vrule width .4pt} \hrule height .4pt}}}

\def\RR{\mathbb{R}}

\def\ZZ{\mathbb{Z}}

\DeclareMathOperator {\var}{var}

\def\clearp{}
\def\xyz#1{#1}

\begin{document}

\title{Spatial heterogeneity can explain the stable\\
 coexistence of savanna and forest in South America }
\author{Ruibo Ma and Rick Durrett\\
Beijing Jiaotong University and Duke University     }

\maketitle

\begin{abstract}
In work with a variety of co-authors, Staver and Levin have argued that savanna and forest coexist as alternative stable states with discontinuous changes in density of trees at the boundary. Here we formulate a nonhomogeneous spatial model of the competition between forest and savanna. We prove that coexistence occurs for a time that is exponential in the size of the system, and that after an initial transient, boundaries between the alternative equilibria remain stable. This result is valid in general for systems that exhibit bistability in  a homogeneous environment.
\end{abstract}

\section{Introduction} \label{sec:introd}

In a paper published in \textit{Science} \cite{SR0}, Carla Staver, Sally Archibald and Simon Levin argued that tree cover does not increase continuously with rainfall but rather is constrained to low (less than $50\%$, ``savanna'') or high (greater than $75\%$, ``forest'') levels. In follow-up work published in \textit{Ecology} \cite{SR1},  \textit{The American Naturalist} \cite{SR2} and \textit{Journal of Mathematical Biology} \cite{SSL}, they studied the following ODE for the evolution of the fraction of land covered by grass $G$, saplings $S$, and trees $T$;
\beq
	\begin{aligned}
	& \frac{dG}{dt}=\mu S+\nu T- \beta GT ,
		 \\
	& \frac{dS}{dt}=\beta GT-\omega(G) S-\mu S ,
	\\
	& \frac{dT}{dt}=\omega(G) S-\nu T. 
	\end{aligned}           \label{mODE}
\eeq
Here, $\mu \ge \nu$ are the death rates for saplings and trees, and $\omega(G)$ is the rate at which saplings grow into trees when the local density of grass is $G$. Fires decrease this rate of progression, and the incidence of fires
is an increasing function of the fraction of grass, so $\omega(G)$ is decreasing.
Studies suggest (see \cite{SR2} for references) that regions with
tree cover below about 40\% burn frequently but fire is rare above this threshold, so they used
an $\omega$ that is close to a step function.

Inspired by this work, Durrett and Zhang \cite{DZ2015} considered two stochastic spatial models in which each site can be in state 0, 1, or 2: (i) Krone's model \cite{Krone} in which 0 = vacant, 1 = juvenile, and 2 = a mature individual capable of giving birth, and (ii) the Staver-Levin foresst model in which 0 = grass,
1 = sapling, and 2 = tree. Theorem 1 in \cite{DZ2015} shows that if $(0,0,1)$ is an unstable fixed point of \eqref{mODE} then when the range of interaction is large,
there is positive probability of grass and trees surviving starting from a finite set and there is a stationary distribution
in which all three types are present. The result they obtain is asymptotically sharp for Krone's model.
However, in the Staver-Levin forest model, if $(0,0,1)$ is attracting then there may also be another stable fixed point for the ODE, 
and as their Theorem 3 shows, in some of these cases there is a nontrivial stationary distribution.

Touboul, Staver, and Levin \cite{TSS} have investigated a number of modifications of the three species system \eqref{mODE}. Variants of
the ODE that add a fourth type called forest trees display a wide variety of behaviors including limit cycles, homoclinic, and heteroclinic orbits.  
Simulations of the spatial version of ODE systems with periodic orbits, see Figure 9 in \cite{Wald}, suggest these systems will have stationary distributions
that are patchy and with local densities oscillating when the scale of observation is smaller than  what physicists call the correlation length (Figure 4 in \cite{DL98}).
Proving that this occurs is a very difficult problem. Here we will instead focus our attention on a two species system analyzed in \cite{TSS}. In our notation it is 
\beq
\frac{dF}{dt} = \phi_1(F) G - \phi_0(G)F\textnormal{, where }  G = 1- F .
\label{TSL}
\eeq
We will later argue that the results we prove hold when the system in a homogeneous enivronment exhibits bistability, so we will use a simple model with this property.
Let $\alpha,\beta >1$ and
\begin{align}
\frac{du}{dt} &= A(1-u) u^\alpha - Bu (1-u)^\beta
\nonumber \\
& = u(1-u) [A u^{\alpha-1} - B(1-u)^{\beta-1}]
.
\label{simpleDE} 
\end{align}
To check bistability note that when $u \to 0$ the term in square brackets converges to
$-B$ so 0 is stable, when $u \to 1$ it converges to $A$ so 1 is stable. Interior fixed points
satisfy
$$ 
\frac{u^{\alpha-1}}{(1-u)^{\beta-1}} = \frac{B}{A}
			.
$$
The right-hand side is strictly increasing from 0 to $\infty$ as $u$ increases from 0 to 1 so there is a unique interior fixed point, which we will call $\nu$. Since 0 and 1 are stable it must be unstable.

To study this system we will introduce a stochastic spatial model. Since our motivation comes from vegetation patterns in South America we will restrict our attention to two dimensions. We will take a limit of the system when the range of interaction $L \to \infty$, so we will formulate the model on a fine lattice ${\cal L} =\ZZ^2/L$. We assume $\alpha$ and $\beta$ are integers in order to have  a simple construction.

\begin{definition} \label{graphrep}
{\bf Graphical representation.}
 Associated with each site $x \in {\cal L}$ there are two Poisson processes.

\begin{itemize}

\item
$B^x_n$, $n \ge 1$ with rate $\lambda$. At time $t= B^x_n$ we have $\alpha$ points $U^{x,1}_n, \ldots U^{x,\alpha}_n$ chosen at random from  $D(x,1) \cap {\cal L} - \{0\}$ without replacement. If $\xi^L_{t-}(x)=0$ and $\xi^L_{t-}(U^{x,i}_n) =1$ for all $1\le i \le \alpha$ then we set $\xi_{t}(x)=1$.

\item
$D^x_n$, $n \ge 1$ with rate $1$. At time $t= D^x_n$ we have $\beta$ points $\bar U^{x,1}_n, \ldots \bar U^{x,\beta}_n$  chosen at random from  $D(x,1) \cap {\cal L} - \{0\}$ without replacement. If $\xi^L_{t-}(x)=1$ and $\xi_{t-}(\bar U^{x,i}_n ) =0$ for all $1\le i \le \beta$ then we set $\xi^L_{t}(x)=0$.

\end{itemize} \end{definition}

\mn
The behavior of the system only depends on the ratio $A/B$ so we will set $B=1$ and $A=\lambda$.

\subsection{Results for the homogeneous system}  \label{sec:resultsh}

To analyze the homogeneous stochastic spatial model, we will use a three-step procedure that has been employed many times. \xyz{We will give references as the story unfolds.}

\mn
\textbf{Step 1.} {\bf Hydrodynamic limit.} {\it Show that the particle system converges to a deterministic limit that is an integro-differential equation (IDE) (or a partial differential equation (PDE))}.

\mn
To state \xyz{our} result, we need to say what it means for the particle system configuration $\xi^L_0(x)$ which consists of only 1's and 0's to converge to the continuous function $u_0(x)$.  To do this in two dimensions, we tile the plane with {\bf small squares} of side length $L^{-\gamma}$ such that the origin is the lower-left corner of one of these squares.

\begin{definition} \label{densdef} 
Let $\gamma>0$ be small. We say that $\xi^L_t$ has {\bf density at least $\sigma$} on the rectangle $R= [a,b] \times [c,d]$, if in each small square $Q \subset R$ the density $\pi(Q) \ge \sigma$.
\end{definition}

\begin{definition} \label{limdef} 
We say that $\xi^L_t \to u(x)$ if for each fixed $K<\infty$ 
$$
\sup_{Q \subset [-K,K]^2, y \in Q} |\pi(Q)- u(y)| \to 0			.
$$
\end{definition}

\begin{theorem} \label{hydrol}
Let $\bar u(t,x)$ be the average value of $u(t,y)$ over $D(x,1)$, the ball of radius 1 around $x$. If $\xi^L_0(x) \to u_0(x)$, a continuous function, then for all $t>0$, we have $\xi^L_t(x) \to u(t, x)$, the solution of
\beq
\frac{du}{dt} =  \lambda(1-u(t,x))\bar u(t,x)^\alpha - u(t,x) (1 - \bar u(t,x))^\beta,
\qquad 
u(0,x) = u_0(x).
\label{homIDE}
\eeq
\end{theorem}

\noindent
The first result of this type \xyz{in the probability literature} was proved by Swindle \cite{GS} for the linear contact process. In this case $\alpha=1$, $\beta=0$. Neuhauser \cite{CNsex} and Besonov and Durrett \cite{BesDur} have studied the quadratic contact process (also called the sexual reproduction model). The version studied in \cite{BesDur} has a limit with $\alpha=2$ and $\beta=0$, while in \cite{CNsex} $(\bar u(t,x))^2$ is replaced by the average of $u^{(2)}(t,y)$, which \xyz{is defined to be} the fraction of adjacent pairs of particles in $D(x,1)$.

\mn
\textbf{Step 2.} \textit{Obtain results that describe the asymptotic behavior of the limiting integro-differential equation.}

\begin{theorem} \label{wspeed}
In one dimension \eqref{homIDE}  has ``traveling wave'' solutions of the form $u(t,x)=w(x-vt)$ with $\lim_{t\to-\infty} w(x)=1$ and  $\lim_{t\to\infty} w(x)=0$. There is a unqiue value $\lambda_0$ so that $v(\lambda_0)=0$. 
\end{theorem}

\mn
The existence and uniqueness \xyz{of travelling wave solutions} is due to Chen \cite{Chen97}. To prove the second claim suppose that $v(\lambda_0)=0$ and we start the system with parameter $\lambda > \lambda_0$ from the stationary wave solution at $\lambda_0$. Standard comparison methods imply that at time 1 the solution dominates a translate of the initial condition so the speed must be positive. \xyz{Q. Do we want to say this is a proof? An alternative is to call this an assumption, or to be cute -- our Riemann Hyposthesis.} 

To obtain results for the asymptotic behavior of the solution we use  Theorems 6.1 and 6.2 in Weinberger  \cite{Hans}. He states his results in $d$-dimensions but here we simplify by only considering $d=2$. To connect with our notation his fixed point $\pi_0=1$, and his unstable fixed point $\pi_0=\rho$.

Intuitively, to understand the asymptotic behavior of the solution, he looks at plane wave solutions of the form $u(x,t) = w_\xi(x\cdot \xi - c^*(\xi)t)$ that move in direction $\xi \in S^1$, the unit circle in the plane. However, to be precise, what he does is to computes the rate of advance of plane wave solutions in a given direction $\xi$, see Section 5 in his paper. In our case since the neighborhood is a ball, all of the speeds are the same $c^*
(\xi)=v$.

To state the theorems from \cite{Hans}, we need to define
$$
{\cal S} = \{ x \in \RR^2: x \cdot \xi \le c^*(\xi) \},
$$
which in general is the limiting shape of the growing solution starting from compactly supported initial distributions.

\mn
\begin{theorem} \label{Wlimit}
Suppose that the interior ${\cal S}^o \neq \emptyset$, i.e., $c^*(\xi)>0$ for all directions $\xi$. Let ${\cal S}^-$ be a closed subset of ${\cal S}^o$. For any $\sigma > \pi_0$ there is a radius $r_\sigma$ so that if $u(0,x) \ge \sigma$ on a ball of radius $r_\sigma$ then 
$$
\liminf_{n\to\infty}  \min_{x \in n {\cal S}^-} u(n,x) =1     .
$$
If $u(0,x)=0$ off of a compact set and ${\cal S}^+$ is a bounded open set containing 
 ${\cal S}$ then 
$$
\liminf_{n\to\infty}  \max_{x \not\in n {\cal S}^+} u(n,x) \le \pi_0          .
$$
\end{theorem}

\mn
As written the last result can be applied to systems with 5 (or more) fixed points. In our situation, where there are only three fixed points, the second conclusion can be improved to show the limit is 0.

\mn
\textbf{Step 3.} {\it Use a block construction to prove results about the stochastic spatial model}. 

\mn
The block construction was developed in 1988 by Bramson and Durrett \cite{BraDur}. It is a method for proving the existence of stationary distributions for particle systems by comparing with mildly dependent oriented percolaiton. An exposition of the method was given in Durrett's 1993 St.~Flour Lecture Notes \cite{StFl}. More examples can be found in Durrett's 2008 Wald Lectures \cite{Wald} and in Cox, Durrett, and Perkins \cite{CDP}, a 2013 paper that made a number of improvements in the methodology. 

In this paper, we use a simple version of the block construction in two spatial dimensions, \xyz{which} is defined by letting $N$ be an integer and $T$ to be a real number that we will specify later. Let $\sigma>\pi_0$ and $r_\sigma$ the constant from Theorem \ref{Wlimit}. In Section 4 we show:

\mn
{\bf Lemma 4.3.}  {\it Suppose $\delta>0$ and $\ep>0$. If we pick the constant $K$ large and then pick $L$ is large enough then whenever $\xi^L_0 $ has density $\ge \sigma$ on $[-N,N]^2$ then with probability $\ge 1-\ep$, $\xi^L_T$ has density $\ge 1-2\delta$ on $[-3N,3N]^2$ when we only use sites whose dual processes do not escape from $[-(K+3)N,(K+3)N]^2$ during $[0,T]$. }

\mn
Combining the last result with a block construction gives:

\begin{theorem} \label{bigdens}
Suppose that the speed $v>0$. For any $\eta>0$ if $L$ is large enough then 
there is a stationary distribution in which the density of 1's is $\ge 1-\eta$. By symmetry if
the speed $v<0$ then for any $\eta>0$ if $L$ is large enough then 
there is a stationary distribution in which the density of 1's is $\le \eta$. 
\end{theorem}

\noindent
By working harder using results in Chapter 5 in Cox, Durrett, and  Perkins \cite{CDP} one can show that if $v>0$ and there are enough 1's in the initial configuration then with high probability the system is $\equiv1$ on a set that grows linearly in time. Combining this with the corrresponding result for $v<0$ then we see that there is no coexistence unless $v=0$.

\subsection{Results for the inhomogeneous system}

We formulate our system on a large torus with a fine lattice $(\ZZ/L \bmod M)^2$. In terms of the application, we assume that the climatic conditions in South American vary slowly on a continental scale. Mathematically we assume that 
\beq
	\lambda(x) = b(x/M) \qquad\hbox{where $b \in {\mathcal  C}^2$}						,
\label{inhomA}
\eeq  
the collection of $f$ that have continuous derivatives  $\partial f/\partial x$, $\partial f/\partial y$, $\partial^2 f/\partial x^2$, $\partial^2 f/\partial x \partial y$ and $\partial ^2 f/\partial y^2$. Our goal is to prove the following claims

\begin{itemize}
  \item In our heterogeneous system it is possible to have stable coexistence of grassland and forest, which, as we have seen, does not occur in a homogeneous system except possibly when $\lambda=\lambda_0$. 
  \item Coexistence occurs because we will have regions that are grassland, i.e., mostly 0's, and regions of forest that contain mostly 1's. 
  \item The boundaries between \xyz{grassland and forest} can be predicted from the birth rate $\lambda(x)$. Specifically, if $\lambda_0$ is the value defined in Theorem \ref{wspeed} where $v(\lambda_0)=0$ then $\{ x : \lambda(x) < 0\}$ is grassland, and $\{ x : \lambda(x) > 0\}$ is forest.  
\end{itemize} 

\noindent
Intuitively, when viewed from an airplane \xyz{(or satellite)} the boundaries between forest and grassland will be stable once equilibrium is reached. However, from the viewpoint of hikers on the ground, there will be transition zones between the two regions. Our methods do not tell \xyz{us} anything about the nature of these transition zones or how their locations change in time.

To prove these conclusions, we take $\lambda_1>\lambda_0$ and look at one connected component $G$ of the open set $\{ x: \lambda(x) >\lambda_1 \} \subset (\RR \bmod M)^2$. To reduce to the homogeneous case we suppose that $\lambda(x) \equiv \lambda_1$ in $G$. Using the block construction results developed in Section 4, we will show that if $L$ is large the forest persists for a long time \xyz{on the part of the fine lattice} in the region $G$ and almost fills it.

\clearp

\section{Proof of Theorem \ref{hydrol}} \label{sec:hydro}

This result is proved in great detail in Neuhauser \cite{CNsex} for the sexual reproduction model and in Cox, Durrett, and Perkins \cite{CDP} for voter model perturbations that converge to PDE, so we will content ourselves to explain the ideas involved, and refer the reader to the two sources just cited for more detail. Recall that in Definition \ref{graphrep} the process was constructed from  a graphical representation. Given this structure, we can define a {\bf dual process} ${\mathcal  I}^{t,x, L}_s$ starting at $x$ at time $t$ and working backwards in time in order to determine the state of the site $x$ at time $t$. Initially ${\mathcal  I}^{t,x, L}_0 = \{ x \}$. Nothing happens until $T_1$, the first time $s$ so that there is a point $y \in \mathcal {I}^{t,x, L}_s$ and an $n$ so that $B^y_n =t-s$ or $D^y_n =t-s$. In the first case, we add $U^{y,1}_n, \ldots, U^{y,\alpha}_n$ to ${\mathcal  I}^{t,x, L}_{T_1}$. In the second case, we add $\bar U^{y,1}_n, \ldots, \bar U^{y,\beta}_n$  to ${\mathcal  I}^{t,x, L}_{T_1}$.  The dual process continues to add points when there is an arrival in a Poisson process at some $y\in \mathcal {I}^{t,x, L}_s$ until we reach time time $t$ in the dual which is time 0 in the original process.

Durrett and Neuhauser \cite{DurNeu94} call $ \mathcal {I}^{t,x, L}_s$ the \textbf{influence set}, since we need to know the states of these sites in order to determine the state of $x$ at time $t$. The analysis of the influence set in this example is particularly simple because points do not move. 

\mn
$\bullet$ If at some time $s \le t$, a point that we want to add to the dual is already in ${\mathcal  I}^{t,x, L}_s$, we say that a \textbf{collision occurs}. This event is denoted by ${\mathcal  C}^{L.x}_t$. 

\mn
$\bullet$ Given two dual processes ${\mathcal  I}^{t,x, L}_s$ and ${\mathcal  I}^{t, y, L}_s$ with $s \le t$, we say that a \textbf{collision occurs between them} if ${\mathcal  I}^{t,x, L}_t \cap {\mathcal  I}^{t, y, L}_t \neq\emptyset$.  Let this event be denoted by ${\mathcal  C}^{L,x,y}_t$. 

\begin{lemma} \label{nohit}
For any fixed $t$, as $L \to\infty$, 
$$
P( {\mathcal  C}^{L,0}_t ), \sup_{x\neq 0} P( {\mathcal  C}^{L,0,x}_t ) 
\le C_1(t)/L^{2/3}\to 0.
$$
\end{lemma}

\begin{proof} By the definitions of ${\mathcal  C}^0_t$ and ${\mathcal  C}^{0,x}_t$, it is enough to bound the probability that a collision occurs in the dual starting from a particle at 0 and another one at $x$. To bound the growth of $|{\mathcal  I}^{t, 0, L}_s|$, we compare with a branching process $Z_t$ with $Z_0=2$ in which particles give birth at rate $\lambda$ to $\alpha$ offspring and at rate $1$ to $\beta$ offspring.  A standard calculation shows that 
\begin{equation}
	E [Z_t] = 2e^{(\lambda\alpha + \beta) t}  \equiv C_1(t).
	\label{meanbd}
\end{equation}
Using Chebyshev's inequality 
$$
P( Z_t \ge L^{2/3}) \le C_1(t)/L^{2/3}.
$$
When $Z_t \le L^{2/3}$ there are $\le L^{2/3}$ births and $\le L^{2/3}$ particles that can be hit by a new born particle. In two dimensions the neighborhood $D(x,1)$ has asymptotically $\pi L^2$ particles so for large $L$ the expected number of collisions is
$\le L^{4/3} / 3 L^2$
which completes the proof. 
\end{proof} 

\begin{proof}[Proof of Theorem \ref{hydrol}] To prove the result we will first prove a version in which the istates of the sites in the initial configuration are independent and have $P(\xi^L_0(x)=1) =u_0(x)$.

\mn
{\it Step 1. Duality implies that $u_L(t,x) = P(\xi^L_t(x)=1)$ converges to a continuous limit $v(t,x)$.}

\mn
As $L\to\infty$ the dual process converges to a branching random walk in which in which particles (i) give birth at rate $\lambda$ to $\alpha$ offspring and at rate $1$ to $\beta$ offspring and the particles born at $x$ are displaced by independent amounts that are uniform on $D(x,1)$. 

This gives us a tree in space-time. We label the births in this tree with $\alpha$ or $\beta$ depending on the type of birth. To compute the state of $x$ at time $t$ we use the 
i.i.d.~random variable that define the initial configuration to assign states to the particles in the branching process at time $t$ in the dual process. We then work up the tree. When we encounter \xyz{an event in one of the two sets of Poisson processes} we change the state as indicated by the dynamics. For example \xyz{suppose, as drawn in Figure 1 the first event we encounter is a birth at $x$. Then if $x$} is in state 0 and all of its children are in state 1 then $x$ flips to state 1. \xyz{If the event is a death at $y$, as in the second event in Figure 1,} then if $y$ is in state 1 and all of its children are in state 0 then $x$ flips to state 0.

\begin{figure}[ht]
\begin{center}
\begin{picture}(260, 240)
\put(90,230){\line(0,-1){200}}
\put(50,160){\line(1,0){80}}
\put(50,160){\line(0,-1){130}}
\put(130,160){\line(0,-1){130}}
\put(130,120){\line(1,0){60}}
\put(90,90){\line(-1,-0){20}}
\put(70,90){\line(0,-1){60}}
\put(190,120){\line(0,-1){90}}
\put(170,60){\line(1,0){50}}
\put(170,60){\line(0,-1){30}}
\put(220,60){\line(0,-1){30}}
\put(30,30){\line(1,0){220}}
\put(48,15){1}
\put(68,15){0}
\put(88,15){1}
\put(128,15){1}
\put(168,15){1}
\put(188,15){0}
\put(218,15){1}
\put(87,158){$\bullet$}
\put(87,88){$\bullet$}
\put(127,118){$\circ$}
\put(187,58){$\bullet$}
\put(95,125){0}
\put(195,90){1}
\put(135,140){1}
\put(95,195){1}
\put(192,47){$x$}
\put(92,75){$y$}
\end{picture}
\caption{A picture of the dual with $\alpha=2$ and $\beta=1$.
Since $\alpha\neq\beta$ we do not have to label the birth events.
Changes of state occur at the black dots but not at the white one.}
\end{center}
\end{figure}

\mn
{\it Step 2. A second moment computation shows that for fixed $t$ the particle system $\xi^L_t(x) \to v(t,x)$
a deterministic limit.}

\mn
We take $\gamma=1/9$, and let $Q_x$ be the small square containing $x$. Let $n=L^{2(1-\gamma)}$ be the number of sites in $Q_x$, let $S_n$ be the number of 1's in $Q_x$ and let $\pi_t(Q_x)$ be the fraction of 1's in $Q_x$.  Using Lemma \ref{nohit} 
$$
\var(S_n) \le n + n^2/L^{2/3}.
$$
Using Chebyshev's inequality
$$
P( |S_n - ES_n| \ge n^{1-\ep}) \le \frac{2n^2/L^{2/3}}{n^{2(1-\ep)} }
= 2n^{2\ep} L^{-2/3}.
$$
If $\gamma = 1/9$, $n=L^{16/9}$ so there are $L^{2/9}$ small squares per unit area.
If we take $\ep=1/16$ then the error bound in the last equation is
$$
L^{4(8/9)/16 -2/3} =L^{-4/9},
$$
so the expected number of errors per unit area is $L^{-2/9} \to 0$. It follows that for any $K< \infty$
$$
\sup_{x \in[-K,K]^2} |\pi_t(Q_x) - v(t,x)) \to 0
$$

\mn
{\it Step 3: A calculation using the infinitesimal generator of the particles system shows that the limit $v(t,x)$ satisfies the IDE.}

\mn
 This step is similar to Section 2.3 of Neuhuaser \cite{CNsex}. Let ${\cal N} = D(0,1) \cap (\ZZ^2/L-\{0\})$
\begin{align}
     \notag
     \frac {d}{ d t} u^L (t, x) 
&=   \lambda (1-u^L(t,x ) ) \left( \frac {1}{| {\cal  N} |} 
\sum_{y \in x  + {\cal  N}  } u^L(t, y) \right)^\alpha \\ 
& \quad   -u^L (t, x)  \left( 1- \frac {1}{|{\cal  N}|} 
\sum_{y \in x + {\cal  N}}   u^L(t, y)\right)^\beta
     \label{IDECalc}
\end{align}
By steps 1 and 2, as $L \to \infty$, the right-hand side converges to 
\begin{equation*}
	\lambda(1-v(t,x))(\bar v(t,x))^\alpha	 - v(t,x)(1 - \bar v(t,x))^\beta
\end{equation*}
where $\bar v(t,y)$ is the average value of $v(t,x)$ over $D(y, 1)$.
\xyz{To interchange sums and limits one needs to use    the Dominated Convergence Theorem, but that is straightforward.}

\mn
{\it Step 4. The conclusion holds when we assume that $\xi^L_0 \to u_0(x)$.}

\mn
This step is simpler when the dual is branching Brownian motion and the limit is a PDE. In that case, we use the random motion of the dual particles in the last $\delta$ units of time before they hit the initial condition to argue that conditional on \xyz{a dual particle landing in a small square then} their distribution is almost uniform over it. See part(d) of the proof in Section 2 of Durrett and Neuhuaser \cite{DurNeu94}

To extend this argument to our situation we construct the dual with jumps on $\RR^2$ \xyz{ and define the particle's location to be the nearest point on $\ZZ^2/L$, ignoring ties since they have zero probability. When we do this, we again have the property that when a dual particle hits the initial condition in a given small square their distribution within that square is almost uniform over it.}  In this way we create an initial condition in which the sates of different small square are equal to that of randomly chosen site within in, a situationto which the previous proof can be applied. 

\xyz{Readers with an eye for problems have undoubtedly noticed the unfortunate fact that uniform jumps and projection onto the lattice do not commute, which means that we do exactly construct the dynamics, but there are only finitely many jumps in the dual out to time $t$, so the difference vanishes in the limit $L \to\infty$.}
\end{proof}

\clearp

\section{Proof of Theorem \ref{Wlimit}} \label{sec:PfWlim}

It suffices to show that Theorems 6.1 and 6.2 in Weinberger \cite{Hans} apply to our system. He studied the asymptotic behavior (as $n\to\infty$) of  discrete time iterations $u_{n+1}(x) = Q [u_n](x)$, where $Q$ acts on functions $u(x): \RR^d \to [0,1]$. Our evolution occurs in continuous time, but our system can be adapted into his setting by defining $Q [u_0] (x)$ as the solution at time 1 of our integro-differential equation \eqref{homIDE}  starting from initial condition $u_0(x)$. \xyz{To be able to treat continuous time examples with his discrete time theory}, Weinberger \cite{Hans} introduced a number of assumptions that are stated on his pages 360--361. The first five \xyz{that we list} are part of his (3.1).

\mn
(i) Let $\mathcal B$ be the collection of continuous functions defined on $\RR^d$ taking values in $[0,\pi_+]$. He assumes $Q[u]\in B$ for all $u \in B$. In practice he only considers $\pi_+=\infty$ or 1 so we suppose $\pi_+=1$ 

\mn
(ii) The operator $Q$ is {\bf translation invariant}. To be precise if $T_yu(x)=u(x+y)$ is the translation operator acting on functions then $Q[T_x[u]] = T_x[Q[u]]$.

\mn
(iii) A consequence of translation invariance is the if $\alpha$ is a constant function then 
$Q[\alpha]$ is again a constant function. We denote the constant value by $g(\alpha)$. He supposes that there are constants $0 \le \pi_0 < \pi_1 \le 1$ so that $g(\alpha)$ is
$$
\begin{cases} < \alpha & \alpha \in (0,\pi_0) \\
> \alpha & \alpha \in (\pi_0,\pi_1) \\ < \alpha & \alpha \in (\pi_1,1] \end{cases}
$$
In our case $\pi_1=1$ and $\pi_0=\rho$ is the unstable fixed point. 

\mn
(iv) $ u \le v$ implies that $Q[u] \le Q[v]$.

\mn
(v) If $u_n \to u$ uniformly on each bounded subset of $\mathbb R^d$, then $Q[u_n](x) \to Q[u](x)$ for all $x$.

\mn
(i)-(iv) are clearly satisifed. (v) follows from the fact that $u(t,x)$ can be written in terms of an expected value of the dual process applied to the intial distribution. See step 1 in Section \ref{sec:hydro}. Since we are deriving a property of the IDE we can suppose that the initial condition is independent sites. 

\medskip
Finally, since we have $\pi_1=1=\pi_+$ Weinberger;s assumption (3.2) is vacuous. \xyz{At this point we have checked that his assumptions are satisfied in our setting, and so we have achieved our goal in this section.}

\clearp

\section{Block construction: homogeneous case} \label{sec:homoblock}

Our proof is based on the argument in Section 4 of Neuhsuaser \cite{CNsex}. She uses $\rho_s$ for the positive stable fixed point and $\rho_u$ for the unstable fixed point. In our case $\rho_s=1$ and $\rho_u=\rho$. The first step is a result about the asymptotic behavior of the limiting intgrodifferential equation. Her Lemma 4.1 in our setting becomes

\begin{lemma} \label{IDE}
Let $\sigma>\rho$, $\delta>0$,  $N > r_\sigma$, and $T=cN$.  If $c$ is chosen large enough and $u(0,x) \ge \sigma$ on $[-N,N]^2$ then $u(T,x) \ge 1 - \delta$ on $[-3N, 3N]$.
\end{lemma}

\mn
Deviating from the order of the steps in  Neuhauser's proof, we next do a spatial truncation. Define the branching random walk $Y_t$ as follows. Let $U^{i}$ be independent and uniform on $D(0,1)$. There are two types of branching events. A particle at $x$ (i) gives birth at rate $\lambda$ to particles at $x+U^{1}, ..., x+U^{\alpha}$ and (ii)  at rate $1$ gives birth to particles at $x+U^{1}, ..., x + U^{\beta}$. To control the movement of the dual, we use the following well-known result (see, e.g., Biggins \cite{B77}) which can be proved easily using large deviations estimates for random walk.

\begin{lemma} \label{BRW}
Let $N$ and $T=cN$ be as in Lemma \ref{IDE}. Let $Y^0_t$ be the branching random walk starting from a single particle at  $0$. If $K$ is sufficiently large, there is an $\eta>0$ such that 
$$
	P( \eta^0_s \subset [-KN,KN]^2   \text{ for all } 0 \leq s \leq T) \geq 1 - \exp(-\eta T)   
$$
\end{lemma} 

\xyz{Due to our use of a fine lattice it is not possible to show that with high probability all of the dual process do not escape from the box.} 
The locations of particles in $[-3N,3N]^2$ at time $T$ whose dual processes escape from the box $[-(K+3)N,(K+3)N]^2$ will be assumed to vacant. If $N$ is large enough this will reduce the density of occupied sites by at most $\delta=\exp(-\eta T)$. Supposing now that the hydrodynamic limit result in Lemma \ref{IDE} holds with $\delta$ replaced by $2\delta$ we let $L \to \infty$ and use the Theorem \ref{hydrol} to conclude.

\begin{lemma} \label{blockcon}
 Let $\ep_K>0$. If $L$ is large enough and $\xi^L_0 $ has density $\ge \sigma$ on $[-N,N]^2$ then with probability $\ge 1-\ep_K$, $\xi^L_T$ has density $\ge 1-2\delta$ on $[-3N,3N]^2$ when we only use sites whose dual processes do not escape from $[-(K+3)N,(K+3)N]^2$ during $[0,T]$. .
\end{lemma}

We are now ready to compare with oriented percolation on ${\cal L}_3 = \{ (\ell,m,n) \in \ZZ^3 : \ell + m +n \hbox{ is even }\}$. We say that $(\ell,m,n)$ is wet if $\xi^L_{nT}$ has density $ \ge 1-2\delta$ on $(2N\ell,2Nm) + [-N,N]^2$ at time $nT$. Theorem 4.3 in Durrett's St.\@ Flour notes \cite{StFl} implies that the wet sites dominate an $(K+3)$-dependent oriented percolation. The result given in Theorem 4.2 for oriented percolation starting from a product measure together standard arguments (see Liggett \cite{Liggett85}) allow us to prove Theorem \ref{bigdens}.

\clearp

\section{Block construction: inhomogeneous case} \label{sec:inhomoblock}

Let $\lambda_1 > \lambda_0$. To flatten out the torus we will
consider $\lambda(x)$ to be a periodic function on $\RR^2$.  Let $G$ be one connected component of the open set $\{ x : \lambda(x) > \lambda_1\}$ and note that due to our definitions the shape of $G$ does not change as $M$ is varied. 
Let $N$ be the constant from Theorem \ref{IDE} which gives the size of the blacks, and let
$$
B_{\ell,m} = (2N\ell,2Nm) + [-N,N]^2
$$
be the blocks used in the block construction. Note that they tile the plane. Let 
$$
I = \{ (\ell,m) : B_{\ell,m} \subset G \}
$$
be the indices of the blocks that fit completely inside of $G$. 

Pick $\ell_1<\ell_2$ and $m_1<m_2$ so that 
$$
J_1 = [\ell_1,\ell_2] \times [m_1,m_2] \subset I
$$
and no larger rectangle has this property. (Even with the assumption of maximality the choice is not unique, \xyz{since the rectangles have different shapes}.) If the size of the torus $M$ is large enough then $J \neq \emptyset$. Mountford \cite{M93} showed

\begin{lemma} \label{surv2d}
There is a constant $\alpha_2$ so that if $M$ is large enough the oriented percolation
on $J$ starting from all sites wet survives for time 
$$
\ge \exp(\alpha_2 (\ell_2-\ell_1) (m_2-m_1) ).
$$
\end{lemma}

\begin{figure}[ht]
\begin{center}
\begin{picture}(210, 180)
\put(30,90){\line(1,-3){20}}
\put(50,30){\line(4,1){80}}
\put(130,50){\line(1,-1){30}}
\put(160,20){\line(1,5){20}}
\put(180,120){\line(-3,1){90}}
\put(90,150){\line(-1,-1){60}}
\put(70,110){$J_1$}
\put(65,125){\line(0,-1){90}}
\put(65,125){\line(1,0){95}}
\put(43,50){\line(1,0){122}}
\put(160,125){\line(0,-1){105}}
\put(75,90){\vector(-1,0){30}}
\put(100,115){\vector(0,1){20}}
\put(150,90){\vector(1,0){20}}
\put(80,60){\vector(0,-1){18}}
\put(150,60){\vector(0,-1){20}}
\end{picture}
\caption{A picture of the construction in the inhomogeneous case.
For convenience $G$ is drawn with a piecewise linear boundary, and we imagine that the
squares for the black construction are very small. Adding the one dimensional pieces in the second step expands $I$ in the directions indicated so that only the small regions in the southwest and southeast corners are not covered.}
\end{center}
\label{fig:inhcon}
\end{figure}

To enlarge the occupied region, we use one dimension block constructions.

\mn
For $\ell_1 \le \ell \le \ell_2$ define $j_1(\ell)$ and $j_2(\ell)$ so that
$$
\{\ell\} \times [j_1(\ell),j_2(\ell)] \subset I.
$$
For $m_1 \le m \le m_2$ define $k_1(m)$ and $k_2(m)$ so that
$$
[k_1(m),k_2(m)] \times \{ m \} \subset I.
$$
Let $J_2$ be the union of the one-dimensional regions just defined. In the example drawn in Figure \ref{fig:inhcon}, $J_2$ contains all of $G$, except for a the small regions near the southwest and southeast corners of $G$.

One can extend the region once again by adding one dimensional chains containing blocks in $J_2-J_1$. However to have \xyz{an argument proof that is simpler to write we will only prove that  with high probability most of $J_2$ is grassland, and we leave it to the reader to contemplate $J_3$, $J_4, \ldots$} 

\xyz{To handle the one-dimensional contact processes involved in the extension from $J_1$ to $J_2$ we use the following result of Durrett and Liu \cite{DL88}. }

\begin{lemma} \label{survd1}
There is a constant $\alpha_2$ so that if $k$ is large enough the oriented percolation
on a one dimensional chain of $k$ blocks starting form all sites occupied survives for time $\ge \exp(\alpha_1 k)$. 
\end{lemma}

\noindent
\xyz{For the proof in this paper} it is useful to know that the survival estimate is valid for the process starting from one site occupied, \xyz{if the contact process survives long enough to process reach} both end point of the interval. On each one dimensional block the rightmost wet site has positive drift and the leftmost wet site has negative drift so in equilibrium the wet sites are within $O(1)$ of the end points. \xyz{For this fact and} much more about oriented percolation, see Durrett \cite{RD84}.

\clearp

\end{document}